# Numerical simulation for fractional Jaulent-Miodek equation associated with energy-dependent Schrodinger potential using two novel techniques

P. Veeresha[1], D. G. Prakasha[1,*], N. Magesh[2], M. M. Nandeppanavar[3], A. John Christopher[2],

**Abstract**

In present work, we investigate the numerical solution of time-fractional Jaulent–Miodek (JM) equations with the aid of two novel techniques namely, coupled fractional reduced differential transform method (CFRDTM) and $q$-homotopy analysis transform method ($q$-HATM). The obtained solutions are presented in a series form, which are converges rapidly. In order to verify the proposed techniques are reliable and accurate, the numerical simulations have been conducted in terms of absolute error. The obtained solutions are presented graphically to ensure the applicability and validity of the considered algorithms. The results of the study reveal that, the $q$-HATM is computationally very effective and accurate as compared to CFRDTM to analyse fractional nonlinear coupled JM equations.

## 1. Introduction

The concept of fractional calculus although 323 years in 2018, has really showed a tremendous consideration and attention in the last 55 years particularly in 1967 Caputo made the first modification under his investigation [1,2]. It has been proved by many pioneers that the integer order models are generalized by fractional order, describe the complex phenomena in a very effective manner [3,4]. The integer order derivatives are local in nature whereas the Caputo fractional derivatives are nonlocal that is, using the integer order derivative, we can analyse the variations in a neighbourhood of a point but by employing the Caputo fractional derivative we can study the changes in complete interval. This fractional differential operator has been applied to study and model the real world problems in many connected areas of science, engineering and technology [5-7]. Moreover, due to the efficient

*Corresponding Author Tel. No.: +91 8095907689
[1]Department of Mathematics, Karnatak University, Dharwad - 580 003, India
[2]P. G. and Research Department of Mathematics, Govt. Arts College for Men, Krishnagiri - 635 001, India
[3]Department of Mathematics, Government College (Autonomous), Sedam Road, Kalaburagi - 585 105, India

**E-mail:** viru0913@gmail.com (P. Veeresha), prakashadg@gmail.com (D. G. Prakasha),
nmagi_2000@yahoo.co.in (N. Magesh), nandeppanavarmm@gmail.com (M. M. Nandeppanavar),
johnjeromemsc1985@gmail.com (A. John Christopher).

and accurate results obtained while applying FC as a mathematical tool in many fields, researcher working in engineering, technology and other branch of science then mathematics has also been trying to point out their viewpoint.

The major merit of FC has proven to be a very suitable tool for describing hereditary memory properties of various phenomena and in additionally, conventional calculus is only a narrow subset of fractional calculus. The mathematical groundwork for derivatives and integrals of fractional order was laid through combined efforts of senior researchers such are Liouville [8], Riemann [9], Caputo [10], Miller and Ross [11], Podlubny [12] and others. The theory of fractional-order calculus has been related to practical projects and it has been applied to chaos theory [13], signal processing [14], electrodynamics [15], finance [16], and other areas. The analytical and numerical solution for differential equations (DE) of arbitrary order arised in the above phenomena plays a critical role in describing the characters of nonlinear problems exist in daily life.

In 1979, Jaulent and Miodek derived the equation as an extension to energy-dependent potentials [17-18] called Jaulent-Miodek equation. The JM equation arised in several associated areas of physics such as fluid mechanics [19], condense matter physics [20], optics [21], and plasma physics [22]. This study is dedicated to investigating the numerical solution for arbitrary order of the following anomalous and challenging physical model,

$$v_t + v_{xxx} + \frac{3}{2}ww_{xxx} + \frac{9}{2}w_x w_{xx} - 6vv_{xx} - 6vww_x - \frac{3}{2}v_x w^2 = 0$$

and

$$w_t + w_{xxx} - 6\, v_x w - 6vw_x - \frac{15}{2}w_x w^2 = 0.$$

(1)

The nonlinear coupled JM equation has been investigated in past years by many authors through distinct techniques like, Sumudu transform homotopy perturbation method (STHPM) [23], Hermite Wavelets and OHAM [24], Invariant subspace method [25], and others [26-28].

On other hand, the last three decades have been witnessed for the discovery of a number of new techniques to explain the non-linear differential system having fractional order, and in parallel to the development of new computational algorithms and symbolic programming. The analytical and numerical solutions for the non-linear fractional differential equations have fundamental importance. Since, most of the complex phenomena are modelled mathematically by non-linear fractional differential equations. In connection with this, Reduced differential transform method (RDTM) is introduced by Keskin and Oturance [29] in 2009, to find solution for linear and nonlinear differential equations [30-31]. The author in [32-33] modified the RDTM to study and find the numerical solution for coupled fractional differential equations, called Coupled fractional differential transform method (CFRDTM). Besides this, in 1992, *Liao Shijun* proposed homotopy analysis method (HAM) [34-35], and it has been effectively employing to find the solution for problems arises in science and technology. The q-HATM was introduced and nurtured by *Singh et al.* [36] and which is an elegant amalgamation of *q*-HAM and Laplace transform.

## 2. Preliminaries

In the present section, we recall some definitions and properties of fractional calculus and Laplace transform:

**Definition 1.** The fractional integral of a function $f(t) \in C_\mu$ ($\mu \geq -1$) and of order $\alpha > 0$, initially defined by Riemann-Liouville and which is presented [8-9] as

$$J^\alpha f(t) = \frac{1}{\Gamma(\alpha)} \int_0^t (t-\tau)^{\alpha-1} f(\tau) d\tau, \quad t > 0, \alpha \in \mathbb{R}^+. \tag{1}$$

**Definition 2.** The fractional derivative of $f \in C_{-1}^n$ in the Caputo [10] sense is defined as

$$D_t^\alpha f(t) = \begin{cases} \frac{d^n f(t)}{dt^n}, & \alpha = n \in \mathbb{N}, \\ \frac{1}{\Gamma(n-\alpha)} \int_0^t (t-\tau)^{n-\alpha-1} f^{(n)}(\tau) d\tau, & n-1 < \alpha < n, n \in \mathbb{N}, \end{cases} \tag{2}$$

where $\alpha$ is the order of the derivative.

**Definition 3.** The Laplace transform ($LT$) of a Caputo fractional derivative $D_t^\alpha f(t)$ is represented [12,37] as

$$L[D_t^\alpha f(t)] = s^\alpha F(s) - \sum_{r=0}^{n-1} s^{\alpha-r-1} f^{(r)}(0^+), \quad (n-1 < \alpha \leq n), \tag{3}$$

where $F(s)$ is symbolize the $LT$ of the function $f(t)$.

**Theorem 1** (Generalized Taylor's formula) [38]

Suppose that $D_a^{k\alpha} f(t) \in C(a,b]$ for $k = 0,1,\ldots,n+1$, where $0 < \alpha \leq 1$, and then we have

$$f(t) = \sum_{i=0}^n \frac{(t-a)^{i\alpha}}{\Gamma[i\alpha+1]} [D_a^{i\alpha} f(t)]_{t=a} + \mathfrak{R}_n^\alpha(t;a), \tag{4}$$

where $\mathfrak{R}_n^\alpha(t;a) = \frac{(t-a)^{(n+1)\alpha}}{\Gamma[(n+1)\alpha+1]} [D_a^{(n+1)\alpha} f(t)]_{t=\xi}$, $a \leq \xi \leq t, \forall t \in (a,b]$ and $D_a^{k\alpha} = D_a^\alpha \cdot D_a^\alpha \cdot D_a^\alpha \cdots D_a^\alpha$ ($k$ times).

## 3. Coupled Fractional Reduced Differential Transform Method (CFRDTM)

In order to present coupled fractional reduced differential transform, $\mathcal{V}(h, k-h)$ is symbolised as the coupled fractional reduced differential transform of $v(x,t)$. If function $v(x, t)$ is analytic and continuously differentiated with respect to time $t$, then we define the fractional coupled reduced differential transform (CFRDT) of $v(x, t)$ as

$$\mathcal{V}(h, k - h) = \frac{1}{\Gamma(h\alpha+(k-h)\beta+1)} \left[D^{(h\alpha+(k-h)\beta)} v(h, k - h)\right]_{t=0}, \quad (5)$$

where as the inverse transform of $\mathcal{V}(h, k - h)$ is presented as

$$v(x, t) = \sum_{k=0}^{\infty} \sum_{h=0}^{k} \mathcal{V}(h, k - h) t^{h\alpha+(k-h)\beta},$$

which is one of the solutions of coupled fractional differential equations.

**Theorem 2 [39]**. Suppose that $U(h, k - h), V(h, k - h)$ and $W(h, k - h)$ are the coupled fractional reduced differential transform of the function $u(x, y), v(x, t)$ and $w(x, t)$, respectively.

 (i) If $u(x, t) = f(x, t) \pm g(x, t)$, then $U(h, k - h) = F(h, k - h) \pm G(h, k - h)$.
 (ii) If $u(x, t) = af(x, t)$, where $a \in \mathbb{R}$, then $U(h, k - h) = aF(h, k - h)$.
 (iii) If $f(x, t) = u(x, t)v(x, y)$, then $F(h, k - h) = \sum_{l=0}^{h} \sum_{s=0}^{k-h} U(h - l, s)V(l, k - h - s)$.
 (iv) If $f(x, t) = D_t^\alpha u(x, t)$, then $F(h, k - h) = \frac{\Gamma((h+1)\alpha+(k-h)\beta+1)}{\Gamma(h\alpha+(k-h)\beta+1)} U(h + 1, k - h)$.
 (v) If $f(x, t) = D_t^\beta v(x, t)$, then $F(h, k - h) = \frac{\Gamma(h\alpha+(k-h+1)\beta+1)}{\Gamma(h\alpha+(k-h)\beta+1)} V(h, k - h + 1)$.

## 4. Fundamental idea of $q$-homotopy analysis transform method ($q$-HATM)

To present the fundamental idea of $q$-HATM [40-43], we consider a fractional order nonlinear partial differential equation of the form:

$$D_t^\alpha \mathcal{V}(x, y, t) + R\,\mathcal{V}(x, y, t) + N\,\mathcal{V}(x, y, t) = f(x, y, t), \quad n - 1 < \alpha \leq n, \quad (6)$$

where $D_t^\alpha \mathcal{V}(x, y, t)$ denotes the Caputo's fractional derivative of the function $\mathcal{V}(x, y, t)$, $R$ is the bounded linear differential operator in $x$ and $t$, (i.e., for a number $\varepsilon > 0$ we have $\|R\mathcal{V}\| \leq \varepsilon \|\mathcal{V}\|$), $N$ specifies the nonlinear differential operator and Lipschitz continuous with $\mu > 0$ satisfying $|N\mathcal{V} - N\mathcal{W}| \leq \mu|\mathcal{V} - \mathcal{W}|$, and $f(x, y, t)$ represents the source term.

Now, by employing the $LT$ on equation (6), we get

$$s^\alpha L[\mathcal{V}(x, y, t)] - \sum_{k=0}^{n-1} s^{\alpha-k-1} \mathcal{V}^{(k)}(x, y, 0) + L[R\mathcal{V}(x, y, t)] + L[N\mathcal{V}(x, y, t)] = L[f(x, y, t)]. \quad (7)$$

On simplifying equation (7), we have

$$L[\mathcal{V}(x, y, t)] - \frac{1}{s^\alpha} \sum_{k=0}^{n-1} s^{\alpha-k-1} \mathcal{V}^k(x, y, 0) + \frac{1}{s^\alpha} \{L[R\mathcal{V}(x, y, t)] + L[N\mathcal{V}(x, y, t)] - L[f(x, y, t)]\} = 0. \quad (8)$$

According to homotopy analysis method [35], the nonlinear operator defined as

$$N[\varphi(x,y,t;q)] = L[\varphi(x,y,t;q)] - \frac{1}{s^\alpha}\sum_{k=0}^{n-1} s^{\alpha-k-1}\varphi^{(k)}(x,y,t;q)(0^+) \qquad (9)$$
$$+ \frac{1}{s^\alpha}\{L[R\varphi(x,y,t;q)] + L[N\varphi(x,y,t;q)] - L[f(x,y,t)]\},$$

where $q \in \left[0,\frac{1}{n}\right]$, and $\varphi(x,y,t;q)$ is real function of $x, y, t$ and $q$.

We construct a homotopy for non-zero auxiliary function as follows:
$$(1 - nq)L[\varphi(x,y,t;q) - \mathcal{V}_0(x,y,t)] = \hbar q N[\varphi(x,y,t;q)], \qquad (10)$$

where $L$ be a symbol of the $LT$, $q \in \left[0,\frac{1}{n}\right]$ ($n \geq 1$) is the embedding parameter, $\hbar \neq 0$ is an auxiliary parameter, $\varphi(x,y,t;q)$ is an unknown function, $\mathcal{V}_0(x,y,t)$ is an initial guess of $\mathcal{V}(x,y,t)$. The following results hold respectively for $q = 0$ and $q = \frac{1}{n}$;
$$\varphi(x,y,t;0) = \mathcal{V}_0(x,y,t), \quad \varphi\left(x,y,t;\frac{1}{n}\right) = \mathcal{V}(x,y,t). \qquad (11)$$

Thus, by amplifying $q$ from 0 to $\frac{1}{n}$, the solution $\varphi(x,y,t;q)$ converge from $\mathcal{V}_0(x,y,t)$ to the solution $\mathcal{V}(x,y,t)$. Expanding the function $\varphi(x,y,t;q)$ in series form by employing Taylor theorem [37] near to $q$, one can get
$$\varphi(x,y,t;q) = \mathcal{V}_0(x,y,t) + \sum_{m=1}^{\infty} \mathcal{V}_m(x,y,t)q^m, \qquad (12)$$

where
$$\mathcal{V}_m(x,y,t) = \frac{1}{m!}\frac{\partial^m \varphi(x,y,t;q)}{\partial q^m}\Big|_{q=0}. \qquad (13)$$

On choosing the auxiliary linear operator, $\mathcal{V}_0(x,y,t)$, $n$ and $\hbar$, the series (12) converges at $q = \frac{1}{n}$ and then it yields one of the solutions for equation (6)
$$\mathcal{V}(x,y,t) = \mathcal{V}_0(x,y,t) + \sum_{m=1}^{\infty} \mathcal{V}_m(x,y,t)\left(\frac{1}{n}\right)^m. \qquad (14)$$

Now, differentiating the *zero-th* order deformation equation (10) $m$-times with respect to $q$ and then dividing by $m!$ and finally taking $q = 0$, which gives
$$L[\mathcal{V}_m(x,y,t) - k_m \mathcal{V}_{m-1}(x,y,t)] = \hbar \mathfrak{R}_m(\vec{\mathcal{V}}_{m-1}), \qquad (15)$$

where
$$\vec{\mathcal{V}}_m = \{\mathcal{V}_0(x,y,t), \mathcal{V}_1(x,y,t), \dots, \mathcal{V}_m(x,y,t)\}. \qquad (16)$$

Employing the inverse $LT$ on equation (15), it yields
$$\mathcal{V}_m(x,y,t) = k_m \mathcal{V}_{m-1}(x,y,t) + \hbar L^{-1}\left[\mathfrak{R}_m(\vec{\mathcal{V}}_{m-1})\right], \qquad (17)$$

where
$$\mathfrak{R}_m(\vec{\mathcal{V}}_{m-1}) = L[\mathcal{V}_{m-1}(x,y,t)] - \left(1 - \frac{k_m}{n}\right)\left(\sum_{k=0}^{n-1} s^{\alpha-k-1}\mathcal{V}^{(k)}(x,y,0) + \frac{1}{s^\alpha}L[f(x,y,t)]\right) \qquad (18)$$
$$+ \frac{1}{s^\alpha}L[R\mathcal{V}_{m-1} + \mathcal{H}_{m-1}]$$

and
$$k_m = \begin{cases} 0, & m \leq 1, \\ n, & m > 1. \end{cases} \qquad (19)$$

In equation (18), $\mathcal{H}_m$ denotes homotopy polynomial and defined as

$$\mathcal{H}_m = \frac{1}{m!}\left[\frac{\partial^m \varphi(x,y,t;q)}{\partial q^m}\right]_{q=0} \text{ and } \varphi(x,y,t;q) = \varphi_0 + q\varphi_1 + q^2\varphi_2 + \cdots. \quad (20)$$

By equations (17) and (18), we have

$$\mathcal{V}_m(x,y,t) = (k_m + \hbar)\mathcal{V}_{m-1}(x,y,t) - \left(1 - \frac{k_m}{n}\right)L^{-1}\left(\sum_{k=0}^{n-1} s^{\alpha-k-1}\mathcal{V}^{(k)}(x,y,0)\right.$$
$$\left. + \frac{1}{s^\alpha}L[f(x,y,t)]\right) + \hbar L^{-1}\left[\frac{1}{s^\alpha}L[R\mathcal{V}_{m-1} + \mathcal{H}_{m-1}]\right]. \quad (21)$$

Finally, on solving equation (21) we get the iterative terms of $\mathcal{V}_m(x,y,t)$. The q-HATM series solution is presented by

$$\mathcal{V}(x,y,t) = \sum_{m=0}^{\infty} \mathcal{V}_m(x,y,t). \quad (22)$$

## 5. CFRDTM solution for time-fractional coupled Jaulent-Miodek equations

Consider the coupled Jaulent-Miodek equations of time-fractional order [23-24]:

$$\frac{\partial^\alpha v}{\partial t^\alpha} + \frac{\partial^3 v}{\partial x^3} + \frac{3}{2}w\frac{\partial^3 w}{\partial x^3} + \frac{9}{2}\frac{\partial w}{\partial x}\frac{\partial^2 w}{\partial x^2} - 6v\frac{\partial v}{\partial x} - 6vw\frac{\partial w}{\partial x} - \frac{3}{2}\frac{\partial v}{\partial x}w^2 = 0, \quad (23)$$

$$\frac{\partial^\beta w}{\partial t^\beta} + \frac{\partial^3 w}{\partial x^3} - 6\frac{\partial v}{\partial x}w - 6v\frac{\partial w}{\partial x} - \frac{15}{2}\frac{\partial w}{\partial x}w^2 = 0, \quad (24)$$

with initial conditions

$$v(x,0) = \frac{1}{8}\lambda^2\left(1 - 4\,\text{sech}^2\left(\frac{\lambda x}{2}\right)\right), \quad (25)$$

$$w(x,0) = \lambda\,\text{sech}\left(\frac{\lambda x}{2}\right), \quad (26)$$

where $\lambda$ is arbitrary constant and $0 < \alpha, \beta \leq 1$. In order to assess the advantages and the accuracy of the CFRDTM, we consider the fractional coupled J-M equations. Firstly, we derive the recursive formula from equations (23) and (24). Now, $V(h, k - h)$ and $W(h, k - h)$ are piked as the CFRDT of $v(x,t)$ and $w(x,t)$, respectively, where $v(x,t)$ and $w(x,t)$ are the solutions of coupled differential equations of fractional order. Here, $V(0,0) = v(x,0)$ and $W(0,0) = w(x,0)$ are the given initial conditions. The following assumptions have been considered without loss of generality:

$$V(0,j) = 0, \quad j = 1,2,3,4,\ldots,$$
$$W(i,0) = 0, \quad i = 1,2,3,4,\ldots,$$

Employing CFRDTM to equation (23), we get the following recursive formula:

$$\frac{\Gamma((h+1)\alpha + (k-h)\beta + 1)}{\Gamma(h\alpha + (k-h)\beta + 1)}V(h+1, k-h)$$
$$= -\left\{\frac{\partial^3}{\partial x^3}V(h, k-h) - \frac{3}{2}\left(\sum_{l=0}^{h}\sum_{s=0}^{k-h} W(h-l,s)\frac{\partial^3}{\partial x^3}W(l, k-h-s)\right)\right. \quad (27)$$

$$+\frac{9}{2}\left(\sum_{l=0}^{h}\sum_{s=0}^{k-h}\frac{\partial}{\partial x}W(h-l,s)\frac{\partial^2}{\partial x^2}W(l,k-h-s)\right) - 6\left(\sum_{l=0}^{h}\sum_{s=0}^{k-h}V(h-l,s)\frac{\partial}{\partial x}V(l,k-h-s)\right)$$

$$-6\left(\sum_{r=0}^{h}\sum_{l=0}^{h-r}\sum_{s=0}^{k-h}\sum_{p=0}^{k-h-s}V(r,k-h-s-p)W(l,s)\frac{\partial}{\partial x}W(h-r-l,p)\right)$$

$$-\frac{3}{2}\left(\sum_{r=0}^{h}\sum_{l=0}^{h-r}\sum_{s=0}^{k-h}\sum_{p=0}^{k-h-s}\frac{\partial}{\partial x}V(r,k-h-s-p)W(l,s)W(h-r-l,p)\right)\Bigg\}.$$

By using the initial condition considered in Equation (25), we get

$$V(0,0) = v(x,0). \qquad (28)$$

In the similar way, we can get the recursive formula from Equation (24) as:

$$\frac{\Gamma(h\alpha + (k-h+1)\beta + 1)}{\Gamma(h\alpha + (k-h)\beta + 1)} W(h, k-h+1) =$$

$$-\Bigg\{\frac{\partial^3}{\partial x^3}W(h,k-h) - 6\left(\sum_{l=0}^{h}\sum_{s=0}^{k-h}\frac{\partial}{\partial x}V(l,k-h-s)W(h-l,s)\right)$$

$$-6\left(\sum_{l=0}^{h}\sum_{s=0}^{k-h}V(l,k-h-s)\frac{\partial}{\partial x}W(h-l,s)\right)$$

$$-\frac{15}{2}\left(\sum_{r=0}^{h}\sum_{l=0}^{h-r}\sum_{s=0}^{k-h}\sum_{p=0}^{k-h-s}\frac{\partial}{\partial x}W(r,k-h-s-p)W(l,s)W(h-r-l,p)\right)\Bigg\}. \qquad (29)$$

In the same manner, from the equation (26), we have

$$W(0,0) = w(x,0). \qquad (30)$$

According to CFRDTM, using recursive equations (27) and (29) with initial conditions respectively considered in equations (28) and (30) simultaneously, we get

$$V(0,0) = \frac{1}{8}\lambda^2\left(1 - 4\,\text{sech}^2\left(\frac{\lambda x}{2}\right)\right),$$

$$W(0,0) = \lambda\,\text{sech}\left(\frac{\lambda x}{2}\right),$$

$$V(1,0) = \frac{2\lambda^5\,\text{csch}^3(\lambda x)\,\text{sech}^4\left(\frac{\lambda x}{2}\right)}{\Gamma[\alpha+1]},$$

$$W(0,1) = \frac{-\lambda^4\,\text{csch}^2(\lambda x)\,\text{sech}^3\left(\frac{\lambda x}{2}\right)}{\Gamma[\beta+1]},$$

$$V(1,1) = -\frac{3\lambda^8(29 - 26\cosh(\lambda x) + \cosh(2\lambda x))\,\text{sech}^6\left(\frac{\lambda x}{2}\right)}{128\,\Gamma[\alpha+\beta+1]},$$

$$W(0,2) = \frac{\lambda^7(-245 + 140\cosh(\lambda x) + \cosh(2\lambda x))\,\text{sech}^5\left(\frac{\lambda x}{2}\right)}{128\,\Gamma[2\beta+1]},$$

$$V(2,0) = \frac{\lambda^8(93 - 74\cosh(\lambda x) + \cosh(2\lambda x))\,\text{sech}^6\left(\frac{\lambda x}{2}\right)}{128\,\Gamma[2\alpha+1]},$$

$W(1,1) = -\frac{3\lambda^7(-5+3\cosh(\lambda x))\operatorname{sech}^5\left(\frac{\lambda x}{2}\right)}{8\Gamma[\alpha+\beta+1]}$,

$\vdots$

The obtained approximate solutions in the series form are presented as

$v(x,t) = \sum_{k=0}^{\infty} \sum_{h=0}^{k} V(h, k-h) t^{h\alpha + (k-h)\beta}$

$= \frac{1}{8}\lambda^2 \left(1 - 4\operatorname{sech}^2\left(\frac{\lambda x}{2}\right)\right) + \frac{2\lambda^5 \operatorname{csch}^3(\lambda x)\operatorname{sech}^4\left(\frac{\lambda x}{2}\right) t^\alpha}{\Gamma[\alpha+1]} + \frac{\lambda^8(93-74\cosh(\lambda x)+\cosh(2\lambda x))\operatorname{sech}^6\left(\frac{\lambda x}{2}\right) t^{2\alpha}}{128\Gamma[2\alpha+1]}$

$- \frac{3\lambda^8(29-26\cosh(\lambda x)+\cosh(2\lambda x))\operatorname{sech}^6\left(\frac{\lambda x}{2}\right) t^{\alpha+\beta}}{128\Gamma[\alpha+\beta+1]} + \cdots$,

$w(x,t) = \sum_{k=0}^{\infty} \sum_{h=0}^{k} W(h, k-h) t^{h\alpha + (k-h)\beta}$

$= \lambda \operatorname{sech}\left(\frac{\lambda x}{2}\right) + \frac{-\lambda^4 \operatorname{csch}^2(\lambda x)\operatorname{sech}^3\left(\frac{\lambda x}{2}\right) t^\beta}{\Gamma[\beta+1]} - \frac{3\lambda^7(-5+3\cosh(\lambda x))\operatorname{sech}^5\left(\frac{\lambda x}{2}\right) t^{\alpha+\beta}}{8\Gamma[\alpha+\beta+1]}$

$+ \frac{\lambda^7(-245+140\cosh(\lambda x)+\cosh(2\lambda x))\operatorname{sech}^5\left(\frac{\lambda x}{2}\right) t^{2\beta}}{128\Gamma[2\beta+1]} + \cdots$.

## 6. $q$-HATM solution for time-fractional coupled Jaulent-Miodek equations

Consider the coupled J-M equations of arbitrary order

$$\begin{cases} \frac{\partial^\alpha v}{\partial t^\alpha} + \frac{\partial^3 v}{\partial x^3} + \frac{3}{2} w \frac{\partial^3 w}{\partial x^3} + \frac{9}{2} \frac{\partial w}{\partial x} \frac{\partial^2 w}{\partial x^2} - 6v \frac{\partial v}{\partial x} - 6vw \frac{\partial w}{\partial x} - \frac{3}{2} \frac{\partial v}{\partial x} w^2 = 0, \\ \frac{\partial^\alpha w}{\partial t^\alpha} + \frac{\partial^3 w}{\partial x^3} - 6 \frac{\partial v}{\partial x} w - 6v \frac{\partial w}{\partial x} - \frac{15}{2} \frac{\partial w}{\partial x} w^2 = 0, \end{cases} \quad 0 < \alpha \leq 1, \tag{31}$$

with initial conditions

$$v(x,0) = \frac{1}{8}\lambda^2 \left(1 - 4\operatorname{sech}^2\left(\frac{\lambda x}{2}\right)\right), \quad w(x,0) = \lambda \operatorname{sech}\left(\frac{\lambda x}{2}\right), \tag{32}$$

Now by performing $LT$ on Equation (31) and make use of conditions provided in Equation (32), we get

$$L[v(x,t)] = \frac{1}{s}\left(\frac{1}{8}\lambda^2 \left(1 - 4\operatorname{sech}^2\left(\frac{\lambda x}{2}\right)\right)\right) + \frac{1}{s^\alpha} L\left\{\frac{\partial^3 v}{\partial x^3} + \frac{3}{2} w \frac{\partial^3 w}{\partial x^3} + \frac{9}{2} \frac{\partial w}{\partial x} \frac{\partial^2 w}{\partial x^2} \right.$$
$$\left. - 6v \frac{\partial v}{\partial x} - 6vw \frac{\partial w}{\partial x} - \frac{3}{2} \frac{\partial v}{\partial x} w^2\right\}, \tag{33}$$

$$L[w(x,t)] = \frac{1}{s}\left(\lambda \operatorname{sech}\left(\frac{\lambda x}{2}\right)\right) + \frac{1}{s^\alpha} L\left\{\frac{\partial^3 w}{\partial x^3} - 6 \frac{\partial v}{\partial x} w - 6v \frac{\partial w}{\partial x} - \frac{15}{2} \frac{\partial w}{\partial x} w^2\right\}.$$

By the help of Equation (33), we define the non-linear operators as

$$N^1[\varphi_1(x,t;q), \varphi_2(x,t;q)] = L[\varphi_1(x,t;q)] - \frac{1}{s}\left(\frac{1}{8}\lambda^2\left(1 - 4\operatorname{sech}^2\left(\frac{\lambda x}{2}\right)\right)\right)$$
$$+ \frac{1}{s^\alpha} L\left\{\frac{\partial^3 \varphi_1(x,t;q)}{\partial x^3} + \frac{3}{2}\varphi_2(x,t;q)\frac{\partial^3 \varphi_2(x,t;q)}{\partial x^3} \right.$$
$$+ \frac{9}{2}\frac{\partial \varphi_2(x,t;q)}{\partial x}\frac{\partial^2 \varphi_2(x,t;q)}{\partial x^2} - 6\varphi_1(x,t;q)\frac{\partial \varphi_1(x,t;q)}{\partial x}$$
$$\left. - 6\varphi_1(x,t;q)\varphi_2(x,t;q)\frac{\partial \varphi_2(x,t;q)}{\partial x} - \frac{3}{2}\frac{\partial \varphi_1(x,t;q)}{\partial x}\varphi_2^2(x,t;q)\right\}, \tag{34}$$

$$N^2[\varphi_1(x,t;q), \varphi_2(x,t;q)] = L[\varphi_2(x,t;q)] - \frac{1}{s}\left(\lambda \operatorname{sech}\left(\frac{\lambda x}{2}\right)\right)$$
$$+ \frac{1}{s^\alpha} L\left\{\frac{\partial^3 \varphi_2(x,t;q)}{\partial x^3} - 6\frac{\partial \varphi_1(x,t;q)}{\partial x}\varphi_2(x,t;q)\right.$$
$$\left. - 6\varphi_1(x,t;q)\frac{\partial \varphi_2(x,t;q)}{\partial x} - \frac{15}{2}\frac{\partial \varphi_2(x,t;q)}{\partial x}\varphi_2^2(x,t;q)\right\}.$$

By applying proposed algorithm, the deformation equation of $m$-$th$ order is given as

$$L[v_m(x,t) - k_m v_{m-1}(x,t)] = \hbar \Re_{1,m}[\vec{v}_{m-1}, \vec{w}_{m-1}],$$
$$L[w_m(x,t) - k_m w_{m-1}(x,t)] = \hbar \Re_{2,m}[\vec{v}_{m-1}, \vec{w}_{m-1}]. \tag{35}$$

where

$$\Re_{1,m}[\vec{v}_{m-1}, \vec{w}_{m-1}] = L[v_{m-1}(x,t)] - \left(1 - \frac{k_m}{n}\right)\frac{1}{s}\left(\frac{1}{8}\lambda^2\left(1 - 4\operatorname{sech}^2\left(\frac{\lambda x}{2}\right)\right)\right)$$
$$+ \frac{1}{s^\alpha} L\left\{\sum_{i=0}^{m-1} v_i \frac{\partial u_{m-1-i}}{\partial x} + \sum_{i=0}^{m-1} v_i \frac{\partial u_{m-1-i}}{\partial x} - \frac{1}{Re}\left(\frac{\partial^2 u_{m-1}}{\partial x^2} + \frac{\partial^2 u_{m-1}}{\partial y^2}\right)\right\},$$
$$\Re_{2,m}[\vec{v}_{m-1}, \vec{w}_{m-1}] = L[w_{m-1}(x,t)] + \left(1 - \frac{k_m}{n}\right)\frac{1}{s}\left(\lambda \operatorname{sech}\left(\frac{\lambda x}{2}\right)\right)$$
$$+ \frac{1}{s^\alpha} L\left\{\sum_{i=0}^{m-1} u_i \frac{\partial v_{m-1-i}}{\partial x} + \sum_{i=0}^{m-1} v_i \frac{\partial v_{m-1-i}}{\partial x} - \frac{1}{Re}\left(\frac{\partial^2 v_{m-1}}{\partial x^2} + \frac{\partial^2 v_{m-1}}{\partial y^2}\right)\right\}. \tag{36}$$

By applying inverse $LT$ on Equation (35), we get

$$v_m(x,t) = k_m v_{m-1}(x,t) + \hbar L^{-1}\{\Re_{1,m}[\vec{v}_{m-1}, \vec{w}_{m-1}]\},$$
$$w_m(x,t) = k_m w_{m-1}(x,t) + \hbar L^{-1}\{\Re_{2,m}[\vec{v}_{m-1}, \vec{w}_{m-1}]\}. \tag{37}$$

On solving above system of equations, we obtain

$$v_0(x,t) = \frac{1}{8}\lambda^2\left(1 - 4\operatorname{sech}^2\left(\frac{\lambda x}{2}\right)\right), w_0(x,t) = \lambda \operatorname{sech}\left(\frac{\lambda x}{2}\right),$$

$$v_1(x,t) = \frac{-2\hbar\lambda^5 \operatorname{csch}^3(\lambda x)\operatorname{sech}^4\left(\frac{\lambda x}{2}\right)t^\alpha}{\Gamma[\alpha+1]}, w_1(x,t) = \frac{\hbar\lambda^4 \operatorname{csch}^2(\lambda x)\operatorname{sech}^3\left(\frac{\lambda x}{2}\right)t^\alpha}{\Gamma[\alpha+1]},$$

$$v_2(x,t) = \frac{-2(n+\hbar)\hbar\lambda^5 \operatorname{csch}^3(\lambda x)\operatorname{sech}^4\left(\frac{\lambda x}{2}\right)t^\alpha}{\Gamma[\alpha+1]} - \frac{\hbar^2\lambda^8(-2+\cosh(\lambda x))\operatorname{sech}^4\left(\frac{\lambda x}{2}\right)t^{2\alpha}}{16\Gamma[2\alpha+1]},$$

$$w_2(x,t) = \frac{(n+\hbar)\hbar\lambda^4 \operatorname{csch}^2(\lambda x)\operatorname{sech}^3\left(\frac{\lambda x}{2}\right)t^\alpha}{\Gamma[\beta+1]} + \frac{\hbar^2\lambda^7(-3+\cosh(\lambda x))\operatorname{sech}^3\left(\frac{\lambda x}{2}\right)t^{2\alpha}}{32\Gamma[2\alpha+1]},$$

$$v_3(x,t) = \frac{-2(n+\hbar)^2\hbar\lambda^5 \operatorname{csch}^3(\lambda x)\operatorname{sech}^4\left(\frac{\lambda x}{2}\right)t^\alpha}{\Gamma[\alpha+1]} - \frac{(n+\hbar)\hbar^2\lambda^8(-2+\cosh(\lambda x))\operatorname{sech}^4\left(\frac{\lambda x}{2}\right)t^{2\alpha}}{16\Gamma[2\alpha+1]}$$
$$- \frac{\hbar^3\lambda^{11}(2(-165+28\cosh(\lambda x)+\cosh(2\lambda x))\Gamma[\alpha+1]^2+3(43-20\cosh(\lambda x)+\cosh(2\lambda x))\Gamma[2\alpha+1])\operatorname{sech}^6\left(\frac{\lambda x}{2}\right)\tanh\left(\frac{\lambda x}{2}\right)t^{3\alpha}}{1024\Gamma[\alpha+1]^2\Gamma[3\alpha+1]},$$

$$w_3(x,t) = \frac{(n+\hbar)^2\hbar\lambda^4 \operatorname{csch}^2(\lambda x)\operatorname{sech}^3\left(\frac{\lambda x}{2}\right)t^\alpha}{\Gamma[\beta+1]} + \frac{(n+\hbar)\hbar^2\lambda^7(-3+\cosh(\lambda x))\operatorname{sech}^3\left(\frac{\lambda x}{2}\right)t^{2\alpha}}{32\Gamma[2\alpha+1]}$$
$$+ \frac{\hbar^3\lambda^{10}((147-92\cosh(\lambda x)+\cosh(2\lambda x))\Gamma[\alpha+1]^2+12(-7+3\cosh(\lambda x))\Gamma[2\alpha+1])\operatorname{sech}^5\left(\frac{\lambda x}{2}\right)\tanh\left(\frac{\lambda x}{2}\right)t^{3\alpha}}{512\,\Gamma[\alpha+1]^2\,\Gamma[3\alpha+1]},$$

$\vdots$

In the similar way, the rest of the iterative terms can be found. Then, the family of $q$-HATM series solution of the of Equation (31) is given by

$$v(x,t) = v_0(x,t) + \sum_{m=1}^{\infty} v_m(x,t)\left(\frac{1}{n}\right)^m,$$
$$w(x,t) = w_0(x,t) + \sum_{m=1}^{\infty} w_m(x,t)\left(\frac{1}{n}\right)^m. \quad (38)$$

If we set $\alpha = 1, \hbar = -1,$ and $n = 1,$ then the obtained solutions $\sum_{m=1}^{N} v_m(x,t)\left(\frac{1}{n}\right)^m$ and $\sum_{m=1}^{N} w_m(x,t)\left(\frac{1}{n}\right)^m$, respectively converges to the exact solutions

$$v(x,t) = \tfrac{1}{8}\lambda^2\left(1 - 4\,\mathrm{sech}^2\left[\tfrac{\lambda}{2}\left(x + \tfrac{1}{2}\lambda^2 t\right)\right]\right) \text{ and } w(x,t) = \lambda\,\mathrm{sech}\left[\tfrac{\lambda}{2}\left(x + \tfrac{1}{2}\lambda^2 t\right)\right]$$

of the classical order coupled Jaulent-Miodek equations as $N \to \infty$.

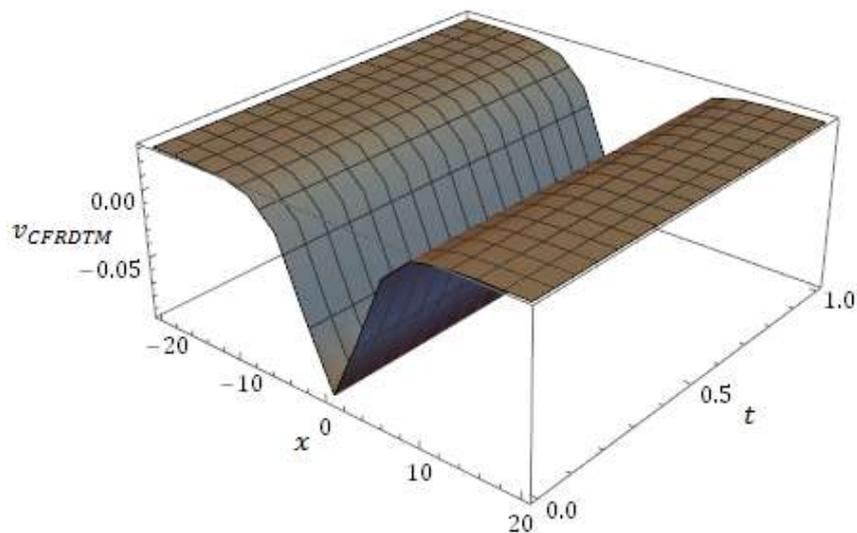

(a)

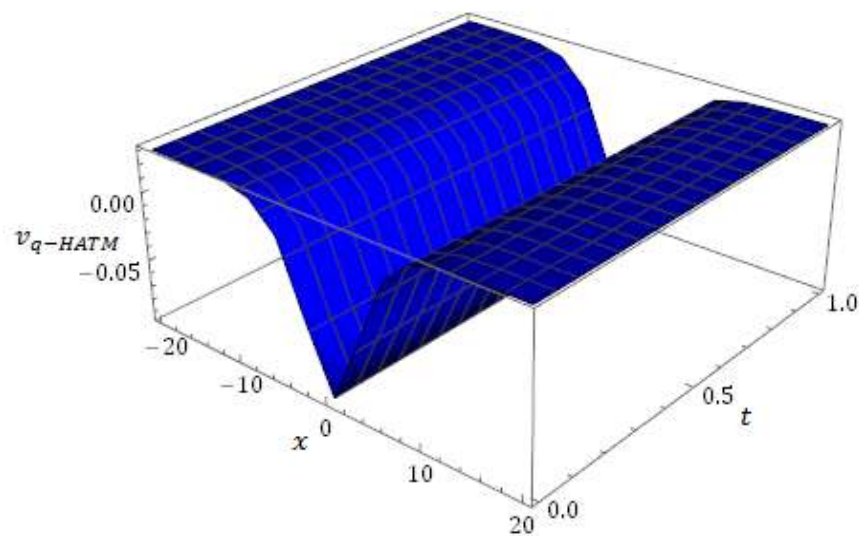

(b)

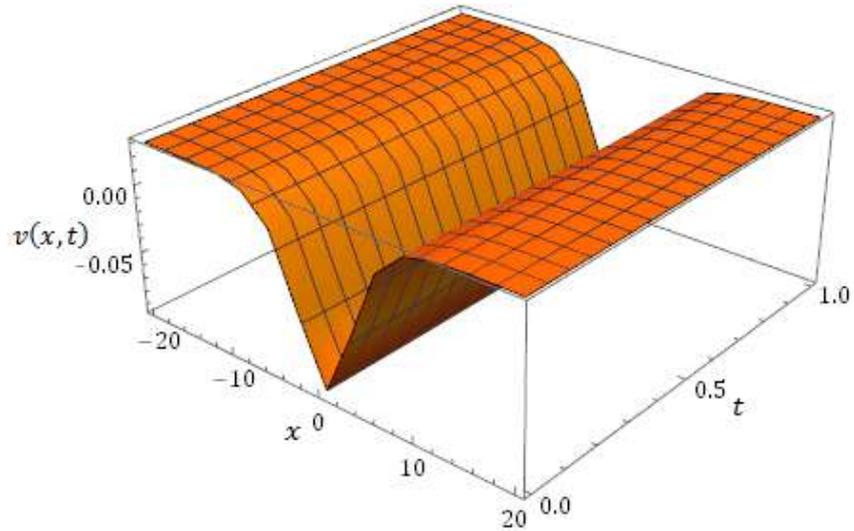

(c)

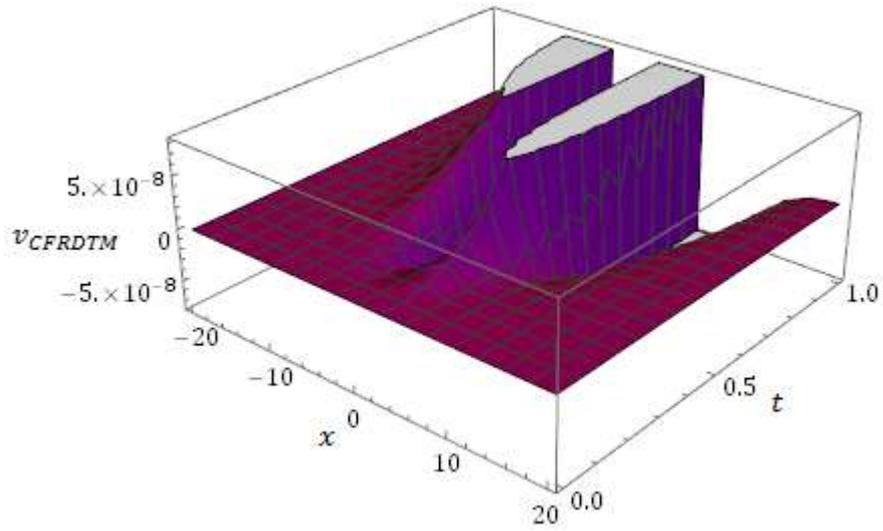

(d)

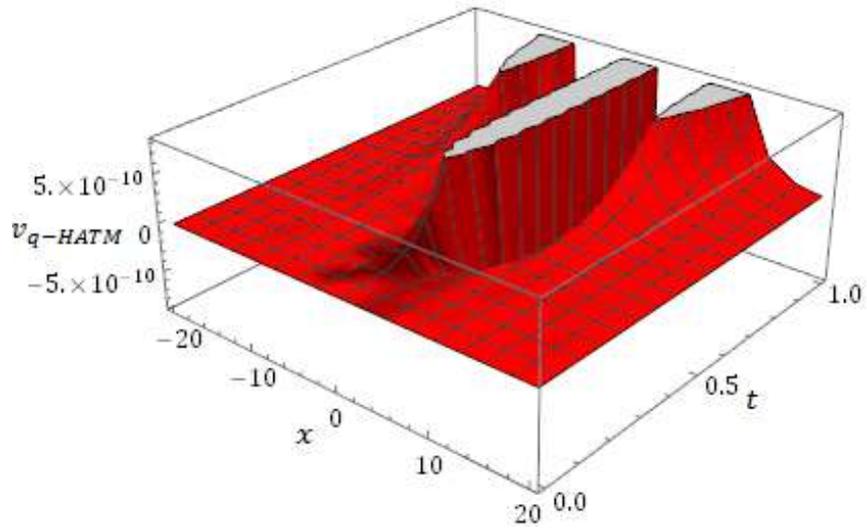

(e)

**Figure 1.** (a) Surface of CFRDTM solution (b) Surface of $q$-HATM solution (c) Surface of Exact solution (d) Surface of Absolute error $= |v_{exa.} - v_{CFRDTM}|$ (e) Surface of Absolute error $= |v_{exa.} - v_{q-HATM}|$ at $\lambda = 0.5, \hbar = -1, n = 1$ and $\alpha = \beta = 1$.

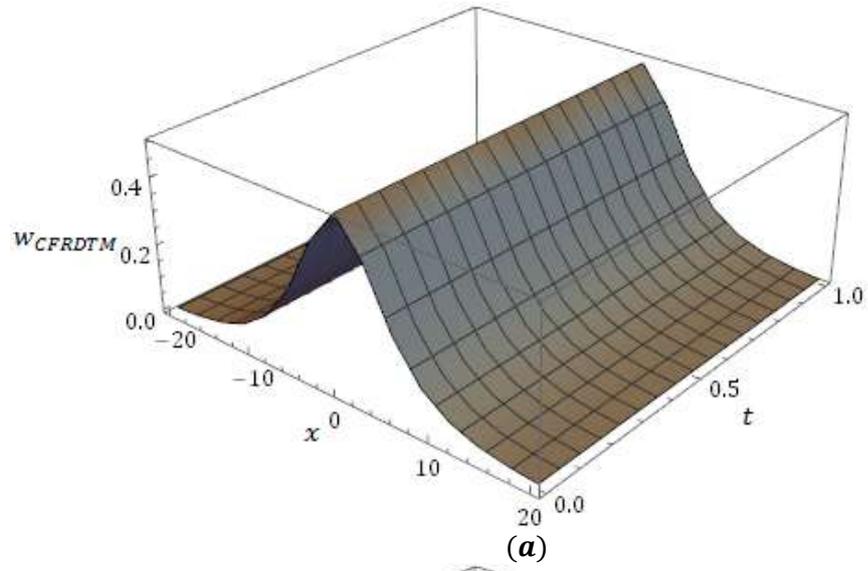

(a)

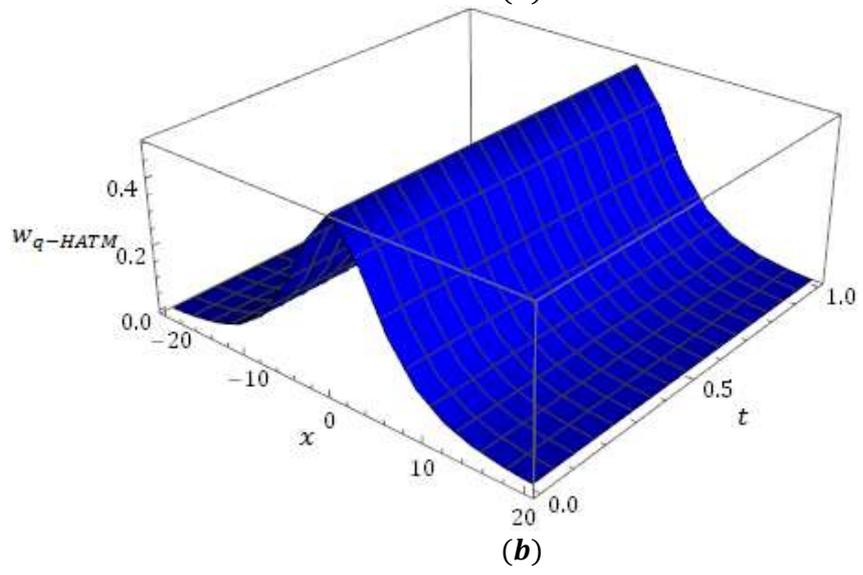

(b)

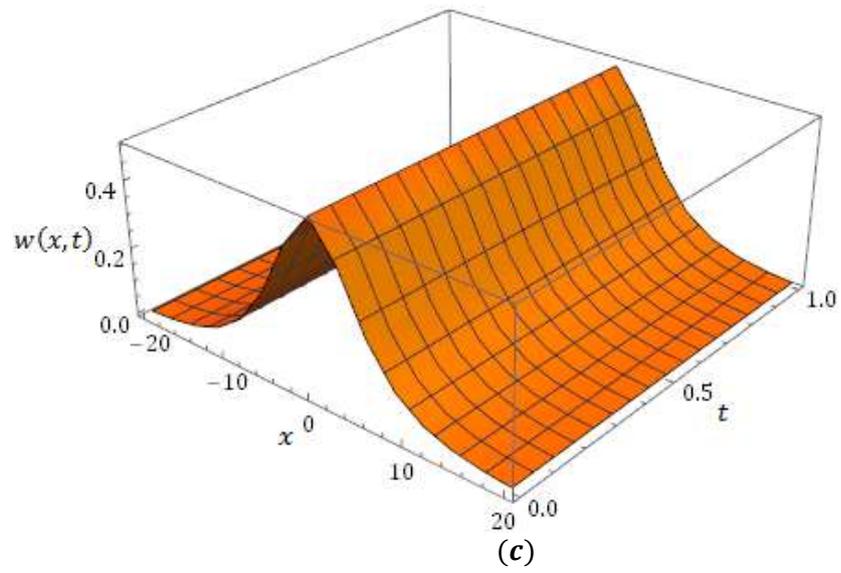

(c)

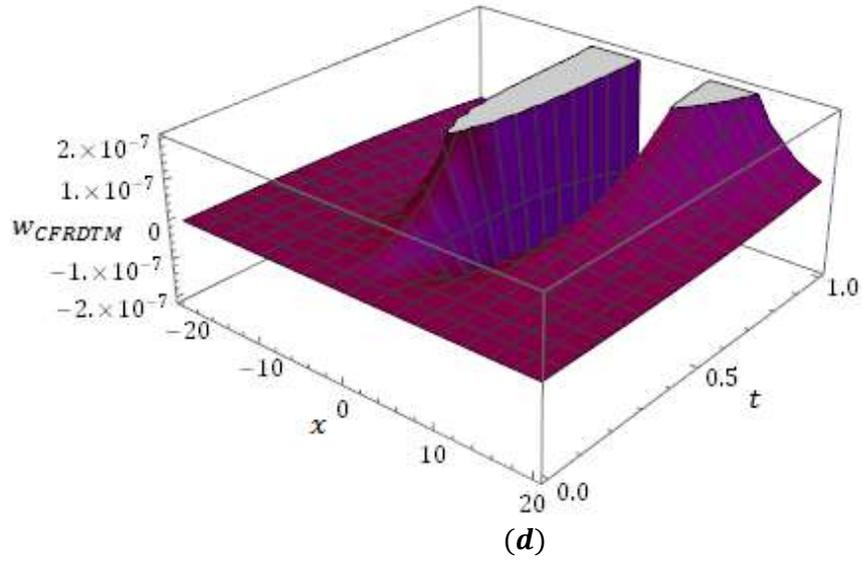

(d)

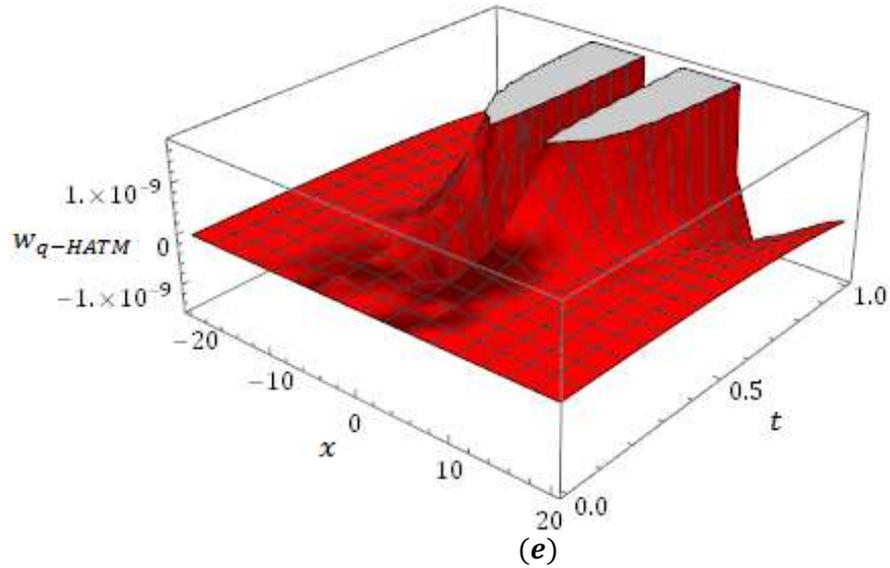

(e)

**Figure 2.** (a) Behaviour of CFRDTM solution (b) Behaviour of $q$-HATM solution (c) Surface of Exact solution (d) Nature of Absolute error $= |w_{exa.} - w_{app.}|$ (e) Nature of Absolute error $= |w_{exa.} - w_{app.}|$ at $\lambda = 0.5, \hbar = -1, n = 1$ and $\alpha = \beta = 1$.

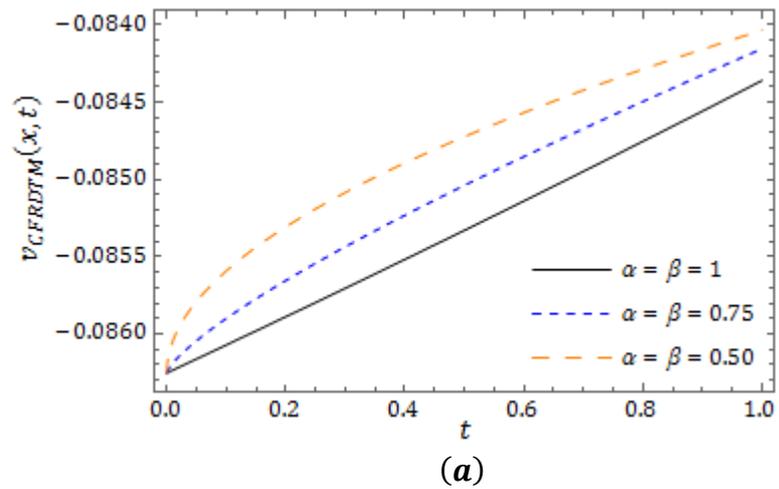

(a)

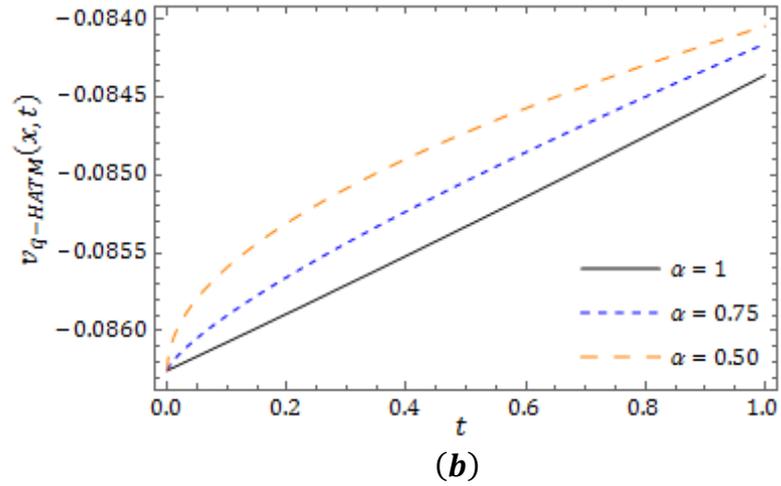

(b)

**Figure 3.** Nature of CFRDTM and $q$-HATM solution $v(x,t)$ with respect to $t$ at $\lambda = 0.5, \hbar = -1, n = 1$ and $x = 1$ with diverse $\alpha$ and $\beta$.

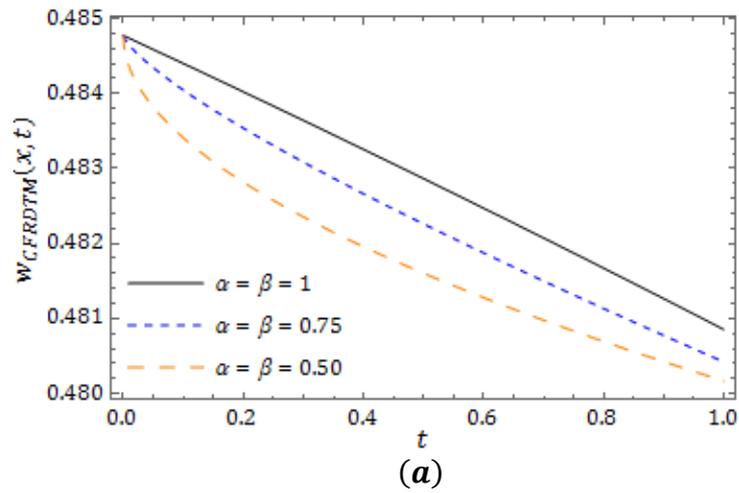

(a)

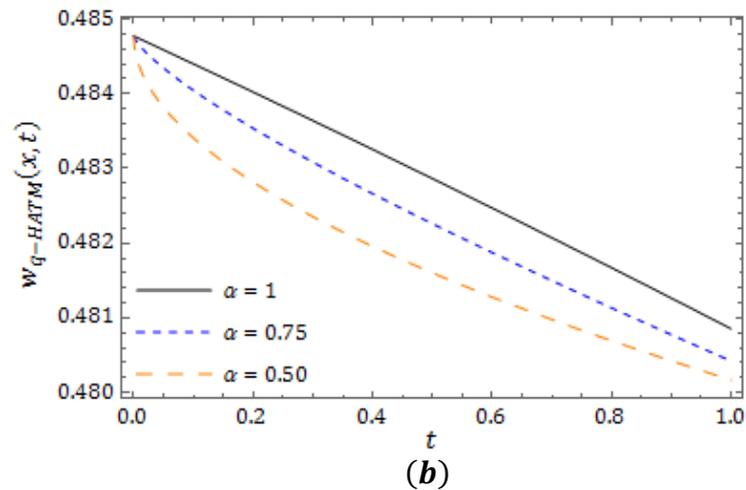

(b)

**Figure 4.** Response of CFRDTM and $q$-HATM solution $w(x,t)$ with respect to $t$ when $\lambda = 0.5, \hbar = -1, n = 1$ and $x = 1$ with diverse $\alpha$ and $\beta$.

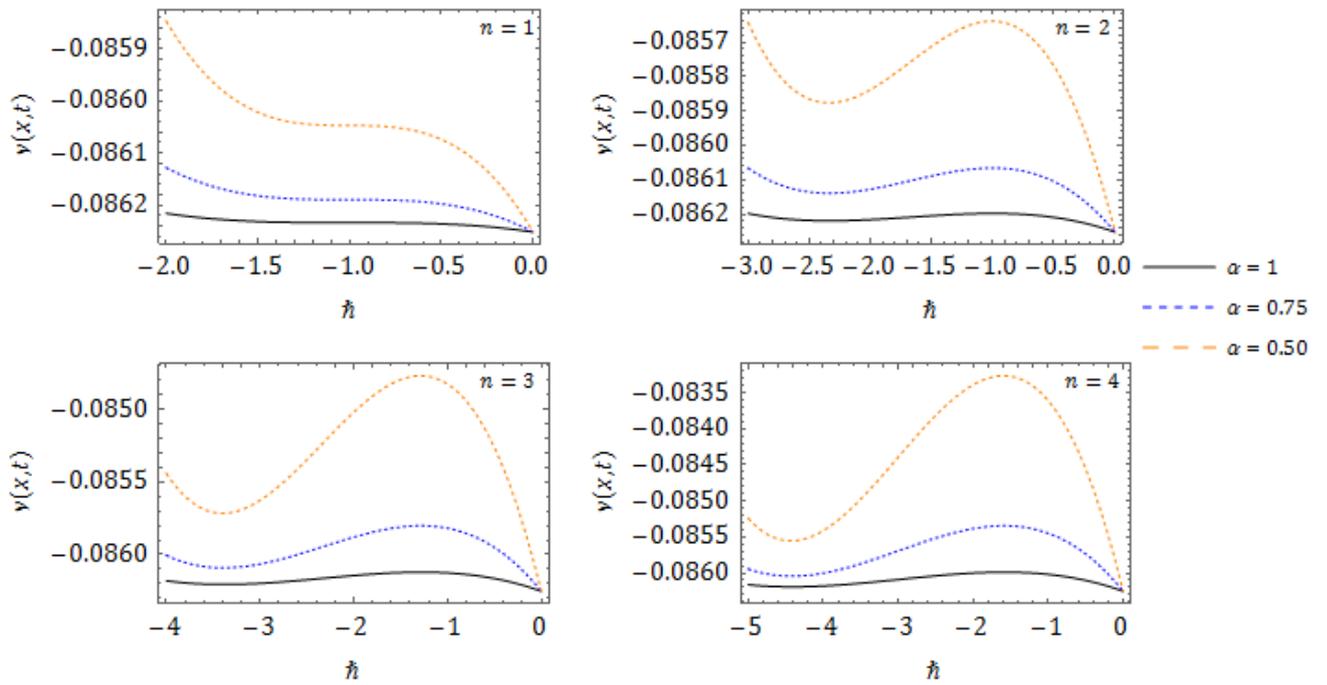

**Figure 5.** $\hbar$-curves drown for $v(x,t)$ with diverse $\alpha$ at $\lambda = 0.5, x = 1$ and $t = 0.01$ for distinct $n$.

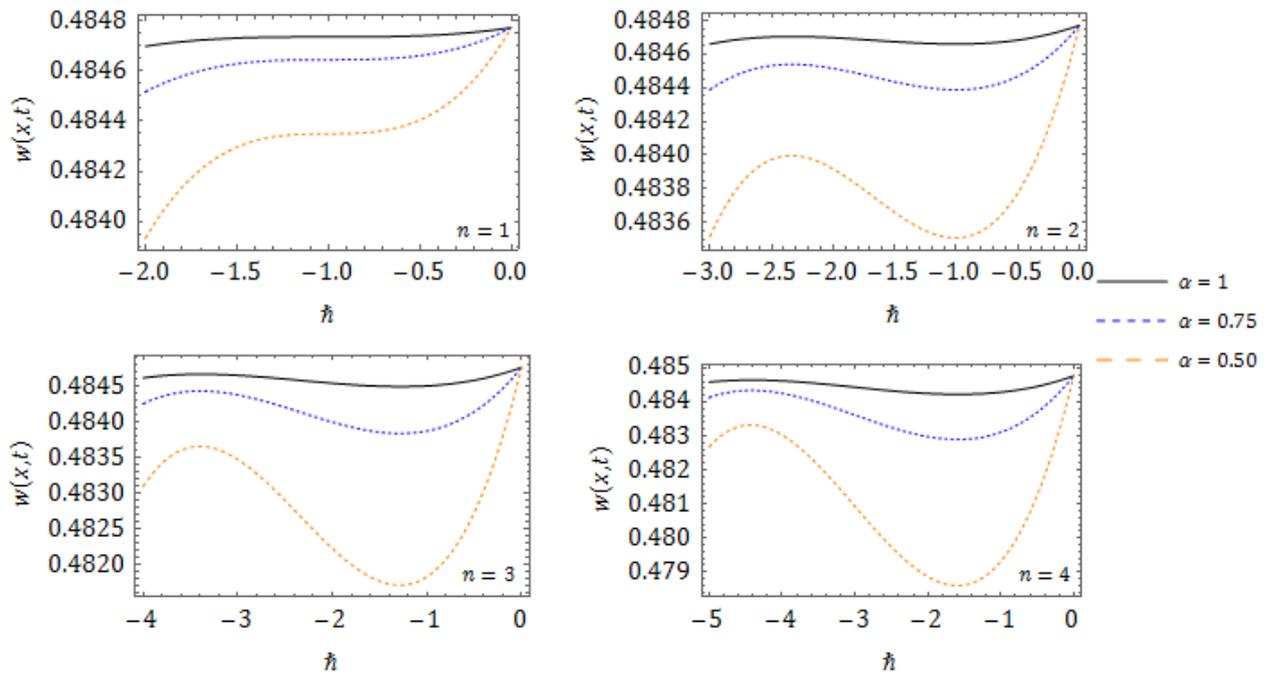

**Figure 6.** $\hbar$-curves drown for $w(x,t)$ with diverse $\alpha$ at $\lambda = 0.5, x = 1$ and $t = 0.01$ for distinct $n$.

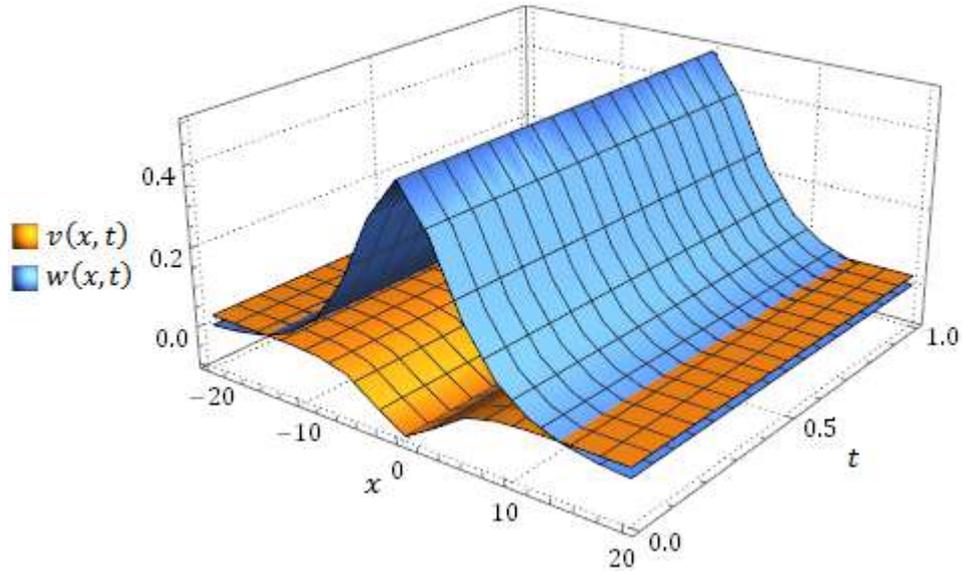

**Figure 7**. Surface of coupled $q$-HATM solutions $v(x,t)$ and $w(x,t)$ at $\lambda = 0.5$, $\hbar = -1$, $n = 1$ and $\alpha = 1$.

**Table 1** Numerical simulation regard to $v(x,t)$ in terms of absolute error obtained by CFRDTM for the Equation (23) with distinct $x$ and $t$ when $\lambda = 0.5$ at diverse value of $\alpha$ and $\beta$.

| $x$ | $t$ | $\alpha = \beta = 0.6$ | $\alpha = \beta = 0.75$ | $\alpha = \beta = 0.9$ | $\alpha = \beta = 1$ |
|---|---|---|---|---|---|
| 0.2 | 0.2 | $1.14998 \times 10^{-4}$ | $6.02573 \times 10^{-5}$ | $2.03596 \times 10^{-5}$ | $4.16463 \times 10^{-9}$ |
|  | 0.4 | $1.49471 \times 10^{-4}$ | $8.40314 \times 10^{-5}$ | $3.01338 \times 10^{-5}$ | $3.43094 \times 10^{-8}$ |
|  | 0.6 | $1.62560 \times 10^{-4}$ | $9.63106 \times 10^{-5}$ | $3.61135 \times 10^{-5}$ | $1.19136 \times 10^{-7}$ |
|  | 0.8 | $1.60446 \times 10^{-4}$ | $9.99188 \times 10^{-5}$ | $3.91164 \times 10^{-5}$ | $2.90296 \times 10^{-7}$ |
|  | 1 | $1.45501 \times 10^{-4}$ | $9.58399 \times 10^{-5}$ | $3.93806 \times 10^{-5}$ | $5.82372 \times 10^{-7}$ |
| 0.4 | 0.2 | $2.00479 \times 10^{-4}$ | $1.07776 \times 10^{-4}$ | $3.71611 \times 10^{-5}$ | $8.02712 \times 10^{-9}$ |
|  | 0.4 | $2.41719 \times 10^{-4}$ | $1.39465 \times 10^{-4}$ | $5.12126 \times 10^{-5}$ | $6.51435 \times 10^{-8}$ |
|  | 0.6 | $2.45798 \times 10^{-4}$ | $1.49292 \times 10^{-4}$ | $5.73798 \times 10^{-5}$ | $2.22972 \times 10^{-7}$ |
|  | 0.8 | $2.26277 \times 10^{-4}$ | $1.44534 \times 10^{-4}$ | $5.80907 \times 10^{-5}$ | $5.35868 \times 10^{-7}$ |
|  | 1 | $1.88523 \times 10^{-4}$ | $1.28109 \times 10^{-4}$ | $5.43349 \times 10^{-5}$ | $1.06088 \times 10^{-6}$ |
| 0.6 | 0.2 | $2.82040 \times 10^{-4}$ | $1.53185 \times 10^{-4}$ | $5.32346 \times 10^{-5}$ | $1.15575 \times 10^{-8}$ |
|  | 0.4 | $3.29252 \times 10^{-4}$ | $1.92175 \times 10^{-4}$ | $7.12894 \times 10^{-5}$ | $9.32845 \times 10^{-8}$ |
|  | 0.6 | $3.24250 \times 10^{-4}$ | $1.99362 \times 10^{-4}$ | $7.75229 \times 10^{-5}$ | $3.17596 \times 10^{-7}$ |
|  | 0.8 | $2.87706 \times 10^{-4}$ | $1.86332 \times 10^{-4}$ | $7.59241 \times 10^{-5}$ | $7.59318 \times 10^{-7}$ |
|  | 1 | $2.27879 \times 10^{-4}$ | $1.57875 \times 10^{-4}$ | $6.82130 \times 10^{-5}$ | $1.49564 \times 10^{-6}$ |
| 0.8 | 0.2 | $3.58201 \times 10^{-4}$ | $1.95656 \times 10^{-4}$ | $6.82860 \times 10^{-5}$ | $1.46221 \times 10^{-8}$ |
|  | 0.4 | $4.10506 \times 10^{-4}$ | $2.41214 \times 10^{-4}$ | $9.00011 \times 10^{-5}$ | $1.17669 \times 10^{-7}$ |
|  | 0.6 | $3.96542 \times 10^{-4}$ | $2.45635 \times 10^{-4}$ | $9.61833 \times 10^{-5}$ | $3.99446 \times 10^{-7}$ |
|  | 0.8 | $3.43691 \times 10^{-4}$ | $2.24590 \times 10^{-4}$ | $9.23029 \times 10^{-5}$ | $9.52257 \times 10^{-7}$ |
|  | 1 | $2.62937 \times 10^{-4}$ | $1.84645 \times 10^{-4}$ | $8.07752 \times 10^{-5}$ | $1.87035 \times 10^{-6}$ |
| 1 | 0.2 | $4.27695 \times 10^{-4}$ | $2.34480 \times 10^{-4}$ | $8.20621 \times 10^{-5}$ | $1.71186 \times 10^{-8}$ |
|  | 0.4 | $4.84169 \times 10^{-4}$ | $2.85780 \times 10^{-4}$ | $1.07039 \times 10^{-4}$ | $1.37489 \times 10^{-7}$ |
|  | 0.6 | $4.61547 \times 10^{-4}$ | $2.87379 \times 10^{-4}$ | $1.13061 \times 10^{-4}$ | $4.65819 \times 10^{-7}$ |
|  | 0.8 | $3.93401 \times 10^{-4}$ | $2.58727 \times 10^{-4}$ | $1.06972 \times 10^{-4}$ | $1.10836 \times 10^{-6}$ |
|  | 1 | $2.93237 \times 10^{-4}$ | $2.08045 \times 10^{-4}$ | $9.18346 \times 10^{-5}$ | $2.17282 \times 10^{-6}$ |

**Table 2** Numerical simulation regard to $w(x,t)$ in terms of absolute error obtained by CFRDTM for the Equation (24) with diverse $x$ and $t$ when $\lambda = 0.5$ at distinct value of $\alpha$ and $\beta$.

| $x$ | $t$ | $\alpha = \beta = 0.6$ | $\alpha = \beta = 0.75$ | $\alpha = \beta = 0.9$ | $\alpha = \beta = 1$ |
|---|---|---|---|---|---|
| 0.2 | 0.2 | $2.30416 \times 10^{-4}$ | $1.20720 \times 10^{-4}$ | $4.07827 \times 10^{-5}$ | $5.21689 \times 10^{-9}$ |
|  | 0.4 | $2.99561 \times 10^{-4}$ | $1.68382 \times 10^{-4}$ | $6.03594 \times 10^{-5}$ | $4.29844 \times 10^{-8}$ |
|  | 0.6 | $3.25817 \times 10^{-4}$ | $1.92979 \times 10^{-4}$ | $7.22983 \times 10^{-5}$ | $1.49281 \times 10^{-7}$ |
|  | 0.8 | $3.21531 \times 10^{-4}$ | $2.00136 \times 10^{-4}$ | $7.82094 \times 10^{-5}$ | $3.63810 \times 10^{-7}$ |
|  | 1 | $2.91426 \times 10^{-4}$ | $1.91790 \times 10^{-4}$ | $7.85409 \times 10^{-5}$ | $7.29980 \times 10^{-7}$ |
| 0.4 | 0.2 | $4.03498 \times 10^{-4}$ | $2.16852 \times 10^{-4}$ | $7.47504 \times 10^{-5}$ | $1.01172 \times 10^{-8}$ |
|  | 0.4 | $4.86882 \times 10^{-4}$ | $2.80810 \times 10^{-4}$ | $1.03056 \times 10^{-4}$ | $8.21274 \times 10^{-8}$ |
|  | 0.6 | $4.95394 \times 10^{-4}$ | $3.00730 \times 10^{-4}$ | $1.15446 \times 10^{-4}$ | $2.81182 \times 10^{-7}$ |
|  | 0.8 | $4.56224 \times 10^{-4}$ | $2.91168 \times 10^{-4}$ | $1.16744 \times 10^{-4}$ | $6.75959 \times 10^{-7}$ |
|  | 1 | $3.80107 \times 10^{-4}$ | $2.57923 \times 10^{-4}$ | $1.08885 \times 10^{-4}$ | $1.33863 \times 10^{-6}$ |
| 0.6 | 0.2 | $5.71624 \times 10^{-4}$ | $3.10320 \times 10^{-4}$ | $1.07799 \times 10^{-4}$ | $1.47159 \times 10^{-8}$ |
|  | 0.4 | $6.68234 \times 10^{-4}$ | $3.89793 \times 10^{-4}$ | $1.44487 \times 10^{-4}$ | $1.18823 \times 10^{-7}$ |
|  | 0.6 | $6.58907 \times 10^{-4}$ | $4.04795 \times 10^{-4}$ | $1.57175 \times 10^{-4}$ | $4.04709 \times 10^{-7}$ |
|  | 0.8 | $5.85337 \times 10^{-4}$ | $3.78632 \times 10^{-4}$ | $1.53839 \times 10^{-4}$ | $9.67991 \times 10^{-7}$ |
|  | 1 | $4.64139 \times 10^{-4}$ | $3.20891 \times 10^{-4}$ | $1.37879 \times 10^{-4}$ | $1.90747 \times 10^{-6}$ |
| 0.8 | 0.2 | $7.32846 \times 10^{-4}$ | $4.00036 \times 10^{-4}$ | $1.39542 \times 10^{-4}$ | $1.88855 \times 10^{-8}$ |
|  | 0.4 | $8.41543 \times 10^{-4}$ | $4.94077 \times 10^{-4}$ | $1.84172 \times 10^{-4}$ | $1.52057 \times 10^{-7}$ |
|  | 0.6 | $8.14513 \times 10^{-4}$ | $5.03993 \times 10^{-4}$ | $1.97007 \times 10^{-4}$ | $5.16457 \times 10^{-7}$ |
|  | 0.8 | $7.07446 \times 10^{-4}$ | $4.61552 \times 10^{-4}$ | $1.89076 \times 10^{-4}$ | $1.23187 \times 10^{-6}$ |
|  | 1 | $5.42633 \times 10^{-4}$ | $3.80010 \times 10^{-4}$ | $1.65197 \times 10^{-4}$ | $2.42089 \times 10^{-6}$ |
| 1 | 0.2 | $8.85417 \times 10^{-4}$ | $4.85021 \times 10^{-4}$ | $1.69634 \times 10^{-4}$ | $2.25214 \times 10^{-8}$ |
|  | 0.4 | $1.00496 \times 10^{-3}$ | $5.92545 \times 10^{-4}$ | $2.21684 \times 10^{-4}$ | $1.80999 \times 10^{-7}$ |
|  | 0.6 | $9.60595 \times 10^{-4}$ | $5.97282 \times 10^{-4}$ | $2.34519 \times 10^{-4}$ | $6.13643 \times 10^{-7}$ |
|  | 0.8 | $8.21325 \times 10^{-4}$ | $5.39082 \times 10^{-4}$ | $2.22087 \times 10^{-4}$ | $1.46106 \times 10^{-6}$ |
|  | 1 | $6.14855 \times 10^{-4}$ | $4.34707 \times 10^{-4}$ | $1.90563 \times 10^{-4}$ | $2.86622 \times 10^{-6}$ |

**Table 3** Numerical study regard to $v(x,t)$ in terms of absolute error obtained by $q$-HATM for system considered in the Equation (31) with diverse $x$ and $t$ when $\lambda = 0.5, \hbar = -1$ and $n = 1,$ at distinct value of $\alpha$.

| $x$ | $t$ | $\alpha = 0.6$ | $\alpha = 0.75$ | $\alpha = 0.9$ | $\alpha = 1$ |
|---|---|---|---|---|---|
| 0.2 | 0.2 | $1.14864 \times 10^{-4}$ | $6.02179 \times 10^{-5}$ | $2.03492 \times 10^{-5}$ | $1.24335 \times 10^{-10}$ |
|  | 0.4 | $1.49003 \times 10^{-4}$ | $8.38441 \times 10^{-5}$ | $3.00661 \times 10^{-5}$ | $1.98709 \times 10^{-9}$ |
|  | 0.6 | $1.61589 \times 10^{-4}$ | $9.58441 \times 10^{-5}$ | $3.59113 \times 10^{-5}$ | $1.00477 \times 10^{-8}$ |
|  | 0.8 | $1.58816 \times 10^{-4}$ | $9.90276 \times 10^{-5}$ | $3.86768 \times 10^{-5}$ | $3.17168 \times 10^{-8}$ |
|  | 1 | $1.43066 \times 10^{-4}$ | $9.43676 \times 10^{-5}$ | $3.85778 \times 10^{-5}$ | $7.73351 \times 10^{-8}$ |
| 0.4 | 0.2 | $2.00217 \times 10^{-4}$ | $1.07699 \times 10^{-4}$ | $3.71407 \times 10^{-5}$ | $1.16347 \times 10^{-10}$ |
|  | 0.4 | $2.40806 \times 10^{-4}$ | $1.39099 \times 10^{-4}$ | $5.10803 \times 10^{-5}$ | $1.85733 \times 10^{-9}$ |
|  | 0.6 | $2.43904 \times 10^{-4}$ | $1.48381 \times 10^{-4}$ | $5.69845 \times 10^{-5}$ | $9.38090 \times 10^{-9}$ |
|  | 0.8 | $2.23097 \times 10^{-4}$ | $1.42794 \times 10^{-4}$ | $5.72312 \times 10^{-5}$ | $2.95781 \times 10^{-8}$ |
|  | 1 | $1.83773 \times 10^{-4}$ | $1.25233 \times 10^{-4}$ | $5.27649 \times 10^{-5}$ | $7.20380 \times 10^{-8}$ |
| 0.6 | 0.2 | $2.81663 \times 10^{-4}$ | $1.53074 \times 10^{-4}$ | $5.32052 \times 10^{-5}$ | $1.03760 \times 10^{-10}$ |
|  | 0.4 | $3.27938 \times 10^{-4}$ | $1.91648 \times 10^{-4}$ | $7.10983 \times 10^{-5}$ | $1.65425 \times 10^{-9}$ |
|  | 0.6 | $3.21523 \times 10^{-4}$ | $1.98048 \times 10^{-4}$ | $7.69518 \times 10^{-5}$ | $8.34447 \times 10^{-9}$ |
| | 0.8 | $2.83130 \times 10^{-4}$ | $1.83823 \times 10^{-4}$ | $7.46824 \times 10^{-5}$ | $2.62762 \times 10^{-8}$ |
| | 1 | $2.21041 \times 10^{-4}$ | $1.53730 \times 10^{-4}$ | $6.59448 \times 10^{-5}$ | $6.39134 \times 10^{-8}$ |
| | 0.2 | $3.57726 \times 10^{-4}$ | $1.95516 \times 10^{-4}$ | $6.82488 \times 10^{-5}$ | $8.74434 \times 10^{-11}$ |
| | 0.4 | $4.08852 \times 10^{-4}$ | $2.40549 \times 10^{-4}$ | $8.97594 \times 10^{-5}$ | $1.39192 \times 10^{-9}$ |
| 0.8 | 0.6 | $3.93110 \times 10^{-4}$ | $2.43979 \times 10^{-4}$ | $9.54607 \times 10^{-5}$ | $7.01008 \times 10^{-9}$ |
| | 0.8 | $3.37931 \times 10^{-4}$ | $2.21425 \times 10^{-4}$ | $9.07319 \times 10^{-5}$ | $2.20392 \times 10^{-8}$ |
| | 1 | $2.54330 \times 10^{-4}$ | $1.79417 \times 10^{-4}$ | $7.79054 \times 10^{-5}$ | $5.35217 \times 10^{-8}$ |
| | 0.2 | $4.27144 \times 10^{-4}$ | $2.34317 \times 10^{-4}$ | $8.20186 \times 10^{-5}$ | $6.84829 \times 10^{-11}$ |
| | 0.4 | $4.82248 \times 10^{-4}$ | $2.85006 \times 10^{-4}$ | $1.06757 \times 10^{-4}$ | $1.08775 \times 10^{-9}$ |
| 1 | 0.6 | $4.57562 \times 10^{-4}$ | $2.85450 \times 10^{-4}$ | $1.12216 \times 10^{-4}$ | $5.46624 \times 10^{-9}$ |
| | 0.8 | $3.86713 \times 10^{-4}$ | $2.55043 \times 10^{-4}$ | $1.05136 \times 10^{-4}$ | $1.71477 \times 10^{-8}$ |
| | 1 | $2.83242 \times 10^{-4}$ | $2.01958 \times 10^{-4}$ | $8.84807 \times 10^{-5}$ | $4.15503 \times 10^{-8}$ |

**Table 4** Numerical study regard to $w(x,t)$ in terms of absolute error obtained by $q$-HATM for system considered in the Equation (31) with diverse $x$ and $t$ when $\lambda = 0.5, \hbar = -1$ and $n = 1$, at distinct value of $\alpha$.

| $x$ | $t$ | $\alpha = 0.6$ | $\alpha = 0.75$ | $\alpha = 0.9$ | $\alpha = 1$ |
|---|---|---|---|---|---|
| | 0.2 | $2.30384 \times 10^{-4}$ | $1.20700 \times 10^{-4}$ | $4.07736 \times 10^{-5}$ | $1.56411 \times 10^{-10}$ |
| | 0.4 | $2.99448 \times 10^{-4}$ | $1.68288 \times 10^{-4}$ | $6.02998 \times 10^{-5}$ | $2.50053 \times 10^{-9}$ |
| 0.2 | 0.6 | $3.25582 \times 10^{-4}$ | $1.92747 \times 10^{-4}$ | $7.21205 \times 10^{-5}$ | $1.26482 \times 10^{-8}$ |
| | 0.8 | $3.21137 \times 10^{-4}$ | $1.99693 \times 10^{-4}$ | $7.78227 \times 10^{-5}$ | $3.99396 \times 10^{-8}$ |
| | 1 | $2.90837 \times 10^{-4}$ | $1.91059 \times 10^{-4}$ | $7.78346 \times 10^{-5}$ | $9.74200 \times 10^{-8}$ |
| | 0.2 | $4.03431 \times 10^{-4}$ | $2.16813 \times 10^{-4}$ | $7.47323 \times 10^{-5}$ | $1.49196 \times 10^{-10}$ |
| | 0.4 | $4.86648 \times 10^{-4}$ | $2.80623 \times 10^{-4}$ | $1.02938 \times 10^{-4}$ | $2.38329 \times 10^{-9}$ |
| 0.4 | 0.6 | $4.94908 \times 10^{-4}$ | $3.00265 \times 10^{-4}$ | $1.15094 \times 10^{-4}$ | $1.20456 \times 10^{-8}$ |
| | 0.8 | $4.55409 \times 10^{-4}$ | $2.90280 \times 10^{-4}$ | $1.15978 \times 10^{-4}$ | $3.80063 \times 10^{-8}$ |
| | 1 | $3.78889 \times 10^{-4}$ | $2.56456 \times 10^{-4}$ | $1.07487 \times 10^{-4}$ | $9.26305 \times 10^{-8}$ |
| | 0.2 | $5.71518 \times 10^{-4}$ | $3.10261 \times 10^{-4}$ | $1.07772 \times 10^{-4}$ | $1.37698 \times 10^{-10}$ |
| | 0.4 | $6.67866 \times 10^{-4}$ | $3.89512 \times 10^{-4}$ | $1.44313 \times 10^{-4}$ | $2.19774 \times 10^{-9}$ |
| 0.6 | 0.6 | $6.58142 \times 10^{-4}$ | $4.04096 \times 10^{-4}$ | $1.56656 \times 10^{-4}$ | $1.10983 \times 10^{-8}$ |
| | 0.8 | $5.84053 \times 10^{-4}$ | $3.77296 \times 10^{-4}$ | $1.52711 \times 10^{-4}$ | $3.49871 \times 10^{-8}$ |
| | 1 | $4.62221 \times 10^{-4}$ | $3.18683 \times 10^{-4}$ | $1.35817 \times 10^{-4}$ | $8.51986 \times 10^{-8}$ |
| | 0.2 | $7.32697 \times 10^{-4}$ | $3.99957 \times 10^{-4}$ | $1.39508 \times 10^{-4}$ | $1.22543 \times 10^{-10}$ |
| | 0.4 | $8.41022 \times 10^{-4}$ | $4.93702 \times 10^{-4}$ | $1.83946 \times 10^{-4}$ | $1.95395 \times 10^{-9}$ |
| 0.8 | 0.6 | $8.13432 \times 10^{-4}$ | $5.03058 \times 10^{-4}$ | $1.96331 \times 10^{-4}$ | $9.85753 \times 10^{-9}$ |
| | 0.8 | $7.05631 \times 10^{-4}$ | $4.59766 \times 10^{-4}$ | $1.87607 \times 10^{-4}$ | $3.10454 \times 10^{-8}$ |
| | 1 | $5.39922 \times 10^{-4}$ | $3.77059 \times 10^{-4}$ | $1.62514 \times 10^{-4}$ | $7.55259 \times 10^{-8}$ |
| | 0.2 | $8.85218 \times 10^{-4}$ | $4.84922 \times 10^{-4}$ | $1.69592 \times 10^{-4}$ | $1.04521 \times 10^{-10}$ |
| | 0.4 | $1.00427 \times 10^{-3}$ | $5.92075 \times 10^{-4}$ | $2.21411 \times 10^{-4}$ | $1.66465 \times 10^{-9}$ |
| 1 | 0.6 | $9.59158 \times 10^{-4}$ | $5.96112 \times 10^{-4}$ | $2.33701 \times 10^{-4}$ | $8.38816 \times 10^{-9}$ |
| | 0.8 | $8.18913 \times 10^{-4}$ | $5.36848 \times 10^{-4}$ | $2.20308 \times 10^{-4}$ | $2.63867 \times 10^{-8}$ |
| | 1 | $6.11250 \times 10^{-4}$ | $4.31016 \times 10^{-4}$ | $1.87313 \times 10^{-4}$ | $6.41164 \times 10^{-8}$ |

**Table 5** Comparative study between CFRDTM [15] and $q$-HATM for $v(x,t)$ and $w(x,t)$ at $\lambda = 0.5, n = 1, \alpha = 1, \beta = 1$ and $\hbar = -1$ with diverse $x$ and $t$.

| $x$ | $t$ | $|v_{Exact} - v_{CFRDTM}|$ | $|v_{Exact} - v_{q-HATM}|$ | $|w_{Exact} - w_{CFRDTM}|$ | $|w_{Exact} - w_{q-HATM}|$ |
|---|---|---|---|---|---|
| 0.2 | 0.1 | $5.12812 \times 10^{-10}$ | $7.77520 \times 10^{-12}$ | $6.42340 \times 10^{-10}$ | $9.77951 \times 10^{-12}$ |
| | 0.2 | $4.16463 \times 10^{-9}$ | $1.24335 \times 10^{-10}$ | $5.21689 \times 10^{-9}$ | $1.56411 \times 10^{-10}$ |
| | 0.3 | $1.42651 \times 10^{-8}$ | $6.29089 \times 10^{-10}$ | $1.78706 \times 10^{-8}$ | $7.91510 \times 10^{-10}$ |
| | 0.4 | $3.43094 \times 10^{-8}$ | $1.98709 \times 10^{-9}$ | $4.29844 \times 10^{-8}$ | $2.50053 \times 10^{-9}$ |
| | 0.5 | $6.79780 \times 10^{-8}$ | $4.84845 \times 10^{-9}$ | $8.51723 \times 10^{-8}$ | $6.10226 \times 10^{-9}$ |
| | 0.6 | $1.19136 \times 10^{-7}$ | $1.00477 \times 10^{-8}$ | $1.49281 \times 10^{-7}$ | $1.26482 \times 10^{-8}$ |
| | 0.7 | $1.91831 \times 10^{-7}$ | $1.86033 \times 10^{-8}$ | $2.40390 \times 10^{-7}$ | $2.34222 \times 10^{-8}$ |
| | 0.8 | $2.90296 \times 10^{-7}$ | $3.17168 \times 10^{-8}$ | $3.63810 \times 10^{-7}$ | $3.99396 \times 10^{-8}$ |
| | 0.9 | $4.18944 \times 10^{-7}$ | $5.07722 \times 10^{-8}$ | $5.25083 \times 10^{-7}$ | $6.39466 \times 10^{-8}$ |
| 0.4 | 0.1 | $9.96126 \times 10^{-10}$ | $7.27987 \times 10^{-12}$ | $1.25533 \times 10^{-9}$ | $9.33215 \times 10^{-12}$ |
| | 0.2 | $8.02712 \times 10^{-9}$ | $1.16347 \times 10^{-10}$ | $1.01172 \times 10^{-8}$ | $1.49196 \times 10^{-10}$ |
| | 0.3 | $2.72872 \times 10^{-8}$ | $5.88342 \times 10^{-10}$ | $3.43967 \times 10^{-8}$ | $7.54698 \times 10^{-10}$ |
| | 0.4 | $6.51435 \times 10^{-8}$ | $1.85733 \times 10^{-9}$ | $8.21274 \times 10^{-8}$ | $2.38329 \times 10^{-9}$ |
| | 0.5 | $1.28135 \times 10^{-7}$ | $4.52925 \times 10^{-9}$ | $1.61564 \times 10^{-7}$ | $5.81382 \times 10^{-9}$ |
| | 0.6 | $2.22972 \times 10^{-7}$ | $9.38090 \times 10^{-9}$ | $2.81182 \times 10^{-7}$ | $1.20456 \times 10^{-8}$ |
| | 0.7 | $3.56533 \times 10^{-7}$ | $1.73588 \times 10^{-8}$ | $4.49676 \times 10^{-7}$ | $2.22974 \times 10^{-8}$ |
| | 0.8 | $5.35868 \times 10^{-7}$ | $2.95781 \times 10^{-8}$ | $6.75959 \times 10^{-7}$ | $3.80063 \times 10^{-8}$ |
| | 0.9 | $7.68191 \times 10^{-7}$ | $4.73215 \times 10^{-8}$ | $9.69162 \times 10^{-7}$ | $6.08270 \times 10^{-8}$ |
| 0.6 | 0.1 | $1.43822 \times 10^{-9}$ | $6.49643 \times 10^{-12}$ | $1.83089 \times 10^{-9}$ | $8.61672 \times 10^{-12}$ |
| | 0.2 | $1.15575 \times 10^{-8}$ | $1.03760 \times 10^{-10}$ | $1.47159 \times 10^{-8}$ | $1.37698 \times 10^{-10}$ |
| | 0.3 | $3.91809 \times 10^{-8}$ | $5.24353 \times 10^{-10}$ | $4.98976 \times 10^{-8}$ | $6.96241 \times 10^{-10}$ |
| | 0.4 | $9.32845 \times 10^{-8}$ | $1.65425 \times 10^{-9}$ | $1.18823 \times 10^{-7}$ | $2.19774 \times 10^{-9}$ |
| | 0.5 | $1.82997 \times 10^{-7}$ | $4.03144 \times 10^{-9}$ | $2.33143 \times 10^{-7}$ | $5.35889 \times 10^{-9}$ |
| | 0.6 | $3.17596 \times 10^{-7}$ | $8.34447 \times 10^{-9}$ | $4.04709 \times 10^{-7}$ | $1.10983 \times 10^{-8}$ |
| | 0.7 | $5.06512 \times 10^{-7}$ | $1.54310 \times 10^{-8}$ | $6.45575 \times 10^{-7}$ | $2.05349 \times 10^{-8}$ |
| | 0.8 | $7.59318 \times 10^{-7}$ | $2.62762 \times 10^{-8}$ | $9.67991 \times 10^{-7}$ | $3.49871 \times 10^{-8}$ |
| | 0.9 | $1.08574 \times 10^{-6}$ | $4.20117 \times 10^{-8}$ | $1.38441 \times 10^{-6}$ | $5.59709 \times 10^{-8}$ |
| 0.8 | 0.1 | $1.82231 \times 10^{-9}$ | $5.47919 \times 10^{-12}$ | $2.35304 \times 10^{-9}$ | $7.67203 \times 10^{-12}$ |
| | 0.2 | $1.46221 \times 10^{-8}$ | $8.74434 \times 10^{-11}$ | $1.88855 \times 10^{-8}$ | $1.22543 \times 10^{-10}$ |
| | 0.3 | $4.94960 \times 10^{-8}$ | $4.41549 \times 10^{-10}$ | $6.39442 \times 10^{-8}$ | $6.19309 \times 10^{-10}$ |
| | 0.4 | $1.17669 \times 10^{-7}$ | $1.39192 \times 10^{-9}$ | $1.52057 \times 10^{-7}$ | $1.95395 \times 10^{-9}$ |
| | 0.5 | $2.30493 \times 10^{-7}$ | $3.38945 \times 10^{-9}$ | $2.97933 \times 10^{-7}$ | $4.76212 \times 10^{-9}$ |
| | 0.6 | $3.99446 \times 10^{-7}$ | $7.01008 \times 10^{-9}$ | $5.16457 \times 10^{-7}$ | $9.85753 \times 10^{-9}$ |
| | 0.7 | $6.36126 \times 10^{-7}$ | $1.29531 \times 10^{-8}$ | $8.22691 \times 10^{-7}$ | $1.82303 \times 10^{-8}$ |
| | 0.8 | $9.52257 \times 10^{-7}$ | $2.20392 \times 10^{-8}$ | $1.23187 \times 10^{-6}$ | $3.10454 \times 10^{-8}$ |
| | 0.9 | $1.35968 \times 10^{-6}$ | $3.52092 \times 10^{-8}$ | $1.75941 \times 10^{-6}$ | $4.96408 \times 10^{-8}$ |
| 1 | 0.1 | $2.13556 \times 10^{-9}$ | $4.29576 \times 10^{-12}$ | $2.80865 \times 10^{-9}$ | $6.54760 \times 10^{-12}$ |
| | 0.2 | $1.71186 \times 10^{-8}$ | $6.84829 \times 10^{-11}$ | $2.25214 \times 10^{-8}$ | $1.04521 \times 10^{-10}$ |
| | 0.3 | $5.78896 \times 10^{-8}$ | $3.45434 \times 10^{-10}$ | $7.61848 \times 10^{-8}$ | $5.27923 \times 10^{-10}$ |
| | 0.4 | $1.37489 \times 10^{-7}$ | $1.08775 \times 10^{-9}$ | $1.80999 \times 10^{-7}$ | $1.66465 \times 10^{-9}$ |
| | 0.5 | $2.69054 \times 10^{-7}$ | $2.64588 \times 10^{-9}$ | $3.54318 \times 10^{-7}$ | $4.05466 \times 10^{-9}$ |
| | 0.6 | $4.65819 \times 10^{-7}$ | $5.46624 \times 10^{-9}$ | $6.13643 \times 10^{-7}$ | $8.38816 \times 10^{-9}$ |
| | 0.7 | $7.41113 \times 10^{-7}$ | $1.00893 \times 10^{-8}$ | $9.76626 \times 10^{-7}$ | $1.55038 \times 10^{-8}$ |
| | 0.8 | $1.10836 \times 10^{-6}$ | $1.71477 \times 10^{-8}$ | $1.46106 \times 10^{-6}$ | $2.63867 \times 10^{-8}$ |
| | 0.9 | $1.58106 \times 10^{-6}$ | $2.73643 \times 10^{-8}$ | $2.08490 \times 10^{-6}$ | $4.21666 \times 10^{-8}$ |

## 7. Numerical results and discussion

In this part, we present numerical study to confirm the accuracy of the numerical solution obtained by the two efficient schemes. In Tables 1 and 2, we present a numerical simulation for fractional coupled JM equations considered in Equation (24) employed by CFRDTM with different $x$ and $t$ for distinct Brownian motion and standard motion ($\alpha = \beta = 1$). Similarly, we conduct a numerical study for coupled system considered in Equation (31) with the aid of an efficient algorithm called $q$-HATM is cited in Tables 3 and 4. Form these tables we can see that, as order increases from 0.6 to 1, the solution get closure to the exact solution. In Table 5, we compared the numerical solutions obtained by two considered computational technique in terms of absolute error.

The behaviour of the solution $v(x,t)$ obtained by CFRDTM and $q$-HATM are presented in Figures 1 (a) and 1 (b) respectively. The nature of the exact solution for JM equation is cited in Figure 1(c) and absolute error for the corresponding equation obtained with the aid of CFRDTM and $q$-HATM are respectively demonstrated in Figures 1 (d) and 1 (e). In a similar manner, the CFRDTM and $q$-HATM solution $w(x,t)$ for Eqs. (24) and (31) is shown in Figure 2. The Figures 3 and 4 are the response of obtained solutions with the help of CFRDTM and $q$-HATM for considered fractional JM equations with distinct Brownian motion and standard motion ($\alpha = \beta = 1$). Figures 5 and 6 demonstrate the $\hbar$-curves for the distinct values of $\alpha$ and $n$ for the coupled system considered in Equation (31) and which help us to adjust and control the convergence region of the series solution. Finally, the coupled surfaces of fractional JM equations consider in Equation (31) are shown in the Figure 7, which guidance us to understand the behaviour of coupled system.

## 8. Conclusion:

In the this study, we have been successfully employed two novel techniques to find the numerical solution of fractional coupled Jaulent–Miodek equations associated with energy-dependent Schrodinger potential. Agreement between numerical results obtained by CFRDTM and $q$-HATM with exact solutions appears very appreciable by means of illustrative results in Tables 1–5. The proposed algorithms are easy to implement, suitable and effective for obtaining the solutions of nonlinear coupled JM equations of fractional order. Moreover, both CFRDTM and $q$-HATM provide the convergent series solutions with easily determinable components without using any perturbation, linearization or limiting assumptions. The numerical simulations illustrate that the results attained with the aid of q-HATM are more exact in comparison with the results achieved by using CFRDTM. Further, $q$-HATM controls and manipulates the series solution, which quickly converges to the exact solution in a short permissible region. Finally, we can conclude the proposed schemes are highly methodical and more accurate, and which can be used to analyse nonlinear problems arisen in complex phenomena.

**Compliance with ethical standards**

Conflict of interest: The authors have declared that no competing interests exist. These authors contributed equally to this work.